\documentclass[11pt]{article}

\usepackage[utf8]{inputenc}
\usepackage{amsmath}
\usepackage{amsfonts,amssymb,latexsym,multicol,color}
\usepackage[francais]{babel}
\usepackage[T1]{fontenc}

\newcounter{fig}
\newtheorem{ques}{Question}

\newtheorem{prop}{Proposition}

\newcommand{\cad}{\text{c'est-\`a-dire }}
\newcommand{\expli}[1]{\quad\text{\footnotesize (#1)}}

\makeatletter
\ifcase\@ptsize
\or
\or
\fi
\makeatother

\newcommand{\implique}{\Rightarrow}

\newcommand{\ioe}{\leqslant}
\newcommand{\soe}{\geqslant}

\newcommand{\vers}{\rightarrow}

\newcommand{\dist}{{\rm dist}}

\newcommand{\Vect}{{\rm Vect}}

\newcommand{\Gram}{{\rm Gram}}

\newcommand{\Ecal}{{\mathcal E}}

\newcommand{\Hcal}{{\mathcal H}}

\newcommand{\Mcal}{{\mathcal M}}

\newcommand{\Scal}{{\mathcal S}}

\providecommand{\scal}[2]{\langle #1, #2 \rangle}

\newcommand{\Nat}{{\mathbb N}}
\newcommand{\Int}{{\mathbb Z}}
\newcommand{\Rat}{{\mathbb Q}}
\newcommand{\Real}{{\mathbb R}}
\newcommand{\Com}{{\mathbb C}}

\newcommand{\virg}{\raisebox{.7mm}{,}}

\newcommand{\fin}{\hfill$\Box$}
\newcommand{\dem}{\noindent {\bf D\'emonstration\ }}
\newcommand{\fine}{\tag*{\mbox{$\Box$}}}

\providecommand{\bysame}{\leavevmode ---\ }
\providecommand{\og}{``}
\providecommand{\fg}{''}
\providecommand{\smfandname}{et}
\providecommand{\smfedsname}{\'eds.}
\providecommand{\smfedname}{\'ed.}
\providecommand{\smfmastersthesisname}{M\'emoire}
\providecommand{\smfphdthesisname}{Th\`ese}

\title{Fonctions arithmétiques multiplicativement monotones}

\author{Michel Balazard}
\date{}

\begin{document}
\maketitle

\begin{center}
  {\sc Abstract}
\end{center}
\begin{quote}
{\footnotesize A real arithmetic function $f$ is \emph{multiplicatively monotonous} if $f(mn)-f(m)$ has constant sign for $m,n$ positive integers. Properties and examples of such functions are discussed, with applications to positive hermitian Toeplitz-multiplicative  determinants.}
\end{quote}

\begin{center}
  {\sc Keywords}
\end{center}
\begin{quote}
{\footnotesize Sets of multiples, Toeplitz-multiplicative determinants, logarithmic density.   \\MSC classification (2010) : 11C20, 11N37 (primary) ; 15A15, 15B05 (secondary).}
\end{quote}

%\newpage

%\tableofcontents

\section{Question initiale}\label{t22}

L'origine du texte qui suit est l'étude de la question suivante.

\begin{ques}\label{t20}
Soit $c : \, \Rat^{+*} \vers \Com$ une fonction telle que
\[
c(1/r) =\overline{c(r)} \quad ( r \in \Rat^{+*}),
\]
et telle que, pour tout $ n \in \Nat$, le déterminant
\[
D_n =\det\big(c(i/j)\big)_{1\ioe i,j \ioe n}
\]
soit strictement positif (par convention, $D_0=1$).

La suite
\[
D_n^{1/n}  \quad ( n \in \Nat^*)
\]
a-t-elle une limite quand $n$ tend vers l'infini?
\end{ques}

\smallskip

Si $c : \Rat^{+*} \vers \Com$ est une fonction quelconque (ne vérifiant pas nécessairement la condition~$c(1/r) =~\overline{c(r)}$), la matrice Toeplitz-multiplicative infinie associée à $c$ est 
\[
M=\big(c(i/j)\big)_{ i,j\soe 1} 
\]

Cette définition, due à Hilberdink (2006, cf. \cite{zbMATH05082236}), est analogue à celle, classique, de matrice de Toeplitz (additive)
\[
\big(c_0(i-j)\big)_{ i,j\soe 1} 
\]
où $c_0$ est définie sur $\Int$ (on numérote alors parfois les lignes et colonnes de la matrice à partir de~$0$, mais nous garderons la même convention pour les deux cas).

\`A la fonction $c$ on peut associer la suite des matrices $M_n$ obtenues à partir de $M$ en ne retenant que les lignes et colonnes d'indices $\ioe n$, où $n \in \Nat^*$. On pose
\[
D_n=\det M_n=\det\big(c(i/j)\big)_{1\ioe i,j\ioe  n} \quad (n \in \Nat^*).
\]
Nous dirons que la suite $(D_n)$ ainsi définie est une suite de déterminants Toeplitz-multiplicatifs.

\smallskip

La condition $c(1/r) =\overline{c(r)}$ assure de plus que les matrices $M_n$ sont hermitiennes, et un critère classique de Sylvester montre que la question posée concerne le cas où les matrices $M_n$ sont toutes définies positives. 

Cette question est motivée par le fait que son analogue pour les déterminants de Toeplitz {\og additifs\fg} a une réponse positive. On a même un résultat plus fort, signalé par Fekete à Szeg\H o, et mentionné dans le premier article de ce dernier (1915, \cite{zbMATH02616579}, \S 3, p. 494-495) : si la fonction $c_0$ vérifie $c_0(-i)=\overline{c_0(i)}$ pour tout $i \in \Int$, et si les déterminants 
\[
\Delta_n=\det\big(c_0(i-j)\big) _{1 \ioe i,j \ioe n}
\]
sont tous strictement positifs, alors la suite $\Delta_n/\Delta_{n-1}$ est décroissante, donc convergente. La démonstration classique de ce fait est rappelée au \S\ref{t3}. Il en résulte que la suite $\Delta_n^{1/n}$ est aussi convergente. Dans le cas des déterminants Toeplitz-multiplicatifs hermitiens positifs, nous verrons au même paragraphe \S\ref{t3} une propriété analogue de décroissance {\og multiplicative\fg} ; avec les hypothèses et les notations de la question ci-dessus, on a
\[
k \mid n \implique D_n/D_{n-1} \ioe D_k/D_{k-1} \quad (k,n \in \Nat^*).
\]

Dès lors, en oubliant provisoirement les spécificités des déterminants Toeplitz-multiplicatifs, il est légitime de s'interroger, dans un cadre plus général, sur les conséquences, pour le comportement asymptotique d'une suite donnée, de l'hypothèse de monotonie multiplicative ; c'est l'objet des paragraphes \S\S\ref{t23}-\ref{t24}. Le résultat principal est l'existence d'une valeur moyenne logarithmique, démontrée à la proposition \ref{t8} au \S\ref{t26}. Après ce détour, je reviendrai aux déterminants Toeplitz-multiplicatifs au \S\ref{t3} ; j'y donnerai un début de réponse à la question initiale, qui restera néanmoins ouverte.

\section{Fonctions arithmétiques multiplicativement monotones}\label{t23}

Au lieu du terme {\og suite\fg}, nous pouvons employer le terme {\og fonction arithmétique\fg}. Il s'agit formellement des mêmes objets ; le changement de vocabulaire induit simplement un changement de point de vue, ici pertinent.

\subsection{Définition}

Nous dirons d'une fonction arithmétique $f$ à valeurs réelles qu'elle est \emph{multiplicativement croissante} si elle vérifie
\[
k \mid n \implique f(k) \ioe f(n) \quad (k,n \in \Nat^*).
\]

La fonction $f$ est \emph{multiplicativement décroissante} si la propriété analogue a lieu avec la dernière inégalité changée de sens ; elle est \emph{multiplicativement monotone} si elle est multiplicativement croissante ou multiplicativement décroissante.

\subsection{Exemples}\label{t4}

Toute fonction arithmétique monotone au sens usuel est multiplicativement monotone. Voici trois autres classes d'exemples de fonctions multiplicativement monotones.

\subsubsection{Ensembles de multiples}

Les fonctions arithmétiques, multiplicativement croissantes, et à valeurs dans $\{0,1\}$, sont exactement les fonctions indicatrices des \emph{ensembles de multiples}. Rappelons que l'ensemble des multiples de la partie $A$ de $\Nat^*$ est défini par
\[
\Mcal(A)=\{ka, \, a\in A, \, k \in \Nat^*\}.
\]

\subsubsection{Fonctions dont la dérivée au sens de Bougaïef est de signe constant}\label{t9}

Bougaïef a défini la {\og dérivée numérique\fg} d'une fonction arithmétique $f$ par la formule 
\[
Df=f*\mu,
\]
(où $*$ désigne la convolution des fonctions arithmétiques, et $\mu$ la fonction de Möbius), et développé la théorie des fonctions arithmétiques à partir de cette notion dans les années 1870 (cf. \cite{bugaev1876} pour un résumé en français). Sa terminologie est imparfaite, notamment parce que cette {\og dérivation\fg} ne vérifie pas de règle de Leibniz\footnote{Il ne faut pas confondre la dérivée au sens de Bougaïef avec la dérivée $f'(n)=f(n)\ln n$ de la théorie des nombres premiers, qui, elle, vérifie la règle de Leibniz $(f*g)'=f'*g+f*g'$.}. Néanmoins, la relation inverse
\[
f(n)=\sum_{d \mid n} (Df)(d),
\]
est décrite par Bougaïef comme le fait que $f$ est l'{\og intégrale numérique\fg} de $Df$ ; ses deux définitions ont donc une certaine cohérence.

Si $Df\soe 0$, alors $f$ est multiplicativement croissante. En effet,
\[
k \mid n \implique f(n)=\sum_{d \mid n} (Df)(d) \soe \sum_{d \mid k} (Df)(d)=f(k).
\]

La réciproque est fausse : la dérivée au sens de Bougaïef de la fonction indicatrice de l'ensemble des multiples de $2$ et $3$ vaut $-1$ pour $n=6$.

\smallskip

\subsubsection{Déterminants Toeplitz-multiplicatifs hermitiens positifs}

Comme nous l'avons signalé au \S \ref{t22}, les quotients $D_n/D_{n-1}$, formés à partir d'une suite de déterminants Toeplitz-multiplicatifs hermitiens positifs, constitue une suite multiplicativement croissante. Cela sera démontré au \S\ref{t3}.

\subsection{Premières propriétés}

Nous énonçons dans ce paragraphe quelques propriétés générales, et immédiates, de la monotonie multiplicative, en nous limitant au cas de la croissance, celui de la décroissance étant analogue. 

\smallskip

{\bf 1.} Une fonction arithmétique $f$ multiplicativement croissante est bornée inférieurement. En effet, on a $f(n)\soe~f(1)$.

{\bf 2.} Si $f$ est une fonction arithmétique multiplicativement croissante et $\alpha$ une fonction à valeurs réelles, définie et croissante sur une partie de $\Real$ contenant l'image de $f$, alors $\alpha\circ f$ est multiplicativement croissante.

{\bf 3.} L'ensemble des fonctions arithmétiques multiplicativement croissantes est un cône convexe.

{\bf 4.} Si $(f_i)_{i \in i}$ est une famille de fonctions arithmétiques multiplicativement croissantes telle que
\[
f(n)=\sup_{i\in I}f_i(n) <\infty \quad (n \in \Nat^*),
\]
alors $f$ est multiplicativement croissante.

{\bf 5.} Si $f$ est une fonction arithmétique à valeurs réelles, il existe une plus petite fonction arithmétique $g$, multiplicativement croissante, et telle que $f \ioe g$. Elle est définie par 
\[
g(n)=\max\{f(d), \, d \mid n\} \quad ( n \in \Nat^*).
\]

Si $f$ est la fonction indicatrice d'une partie $A$ de $\Nat^*$, alors $g$ est la fonction indicatrice de l'ensemble de multiples $\Mcal(A)$.

{\bf 6.} Si $f$ et $g$ sont deux fonctions arithmétiques à valeurs positives ou nulles, et si l'une d'elles, disons $g$, est multiplicativement croissante, alors $f*g$ est multiplicativement croissante. Cela généralise l'exemple du \S\ref{t9} :
\[
k \mid n \implique f*g(n)=\sum_{d \mid n} f(d)g(n/d) \soe \sum_{d \mid k} f(d)g(n/d) \soe \sum_{d \mid k} f(d)g(k/d)=f*g(k).
\]

\section{Facteurs directs de l'ensemble des nombres entiers positifs}

Indépendamment de l'hypothèse de monotonie multiplicative, nous rappelons dans ce paragraphe un cas important d'existence de la limite des moyennes de Cesàro d'une fonction arithmétique. Cette condition suffisante s'appuie sur la notion de \emph{facteur direct}, et sera employée dans la démonstration des propositions \ref{t8} et \ref{t30}.

Une partie $A$ de l'ensemble $\Nat^*$ des nombres entiers positifs est un facteur direct s'il existe une partie $B \subset \Nat^*$ telle que tout nombre entier $n \soe 1$ admette une décomposition unique,
\[
n=a\cdot b \quad (a \in A, \, b \in B).
\]
La partie $B$, dite \emph{complémentaire} de $A$, est manifestement aussi un facteur direct, et $A$ est complémentaire de $B$.

Un théorème d'Erd\H os, Saffari et Vaughan (cf. \cite{zbMATH03520506}, \cite{zbMATH03636131} et \cite{zbMATH03636132}) affirme que tout facteur direct admet une densité asymptotique, égale à la somme des inverses des éléments de son complémentaire :
\begin{equation}\label{t32}
\frac 1x\sum_{\substack{b \in B\\b \ioe x}} 1 \vers \Big(\sum_{a\in A}\frac 1a\Big)^{-1} \quad ( x \vers \infty)
\end{equation}
(si la série diverge, la densité est nulle).

\smallskip

Si $f$ est une fonction arithmétique et $A$, $B$ des facteurs directs complémentaires, la fonction arithmétique $A$-réduite de $f$ est définie par
\begin{equation}\label{t34}
f(n;A)= f(a) \quad ( n=ab, \, a \in A, \, b \in B).
\end{equation}

\begin{prop}\label{t31}
Soit $A$ et $B$ deux facteurs directs complémentaires de $\Nat^*$. On suppose que la somme des inverses des éléments de $A$ converge, et on pose
\[
\lambda = \Big(\sum_{a\in A}\frac 1a\Big)^{-1} (\in\, ]0,1]).
\]

Soit $f$ une fonction arithmétique à valeurs réelles, telle que
\begin{equation}\label{t35}
\sum_{a \in A} \frac{\min\big(f(a),0\big)}a > -\infty.
\end{equation}

Alors
\[
t^{-1}\sum_{n\ioe t} f(n;A) \vers \alpha(f;A) \quad ( t \vers \infty),
\]
où
\[
\alpha(f;A)=\lambda \sum_{a \in A}\frac{f(a)}a \quad (\in\, ]-\infty,\infty]).
\]
\end{prop}
\dem

En notant $a$ un élément générique de $A$, et $b$ un élément générique de $B$, on a
\begin{align*}
t^{-1}\sum_{n\ioe t} f(n;A) &= t^{-1}\sum_{ab\ioe t} f(a)\\
&=\sum_a\frac{f(a)}a (t/a)^{-1}\sum_{b\ioe t/a}1.
\end{align*}

Si $\sum_{a}\lvert f(a)\rvert /a<\infty$, le résultat découle donc de \eqref{t32} et du théorème de Tannery (théorème de convergence dominée de Lebesgue pour les séries). Sinon, il découle de \eqref{t35} et du lemme de Fatou (ou d'un raisonnement élémentaire).\fin

\medskip

Par sommation partielle, on obtient aussi la relation
\begin{equation*}
(\ln t)^{-1}\sum_{n\ioe t} \frac{f(n;A)}n \vers \alpha(f;A) \quad ( t \vers \infty),
\end{equation*}
si $f$ vérifie \eqref{t35}.

\section{\'Etude de la moyenne logarithmique d'une fonction multiplicativement monotone}

La proposition \ref{t8} ci-dessous généralise aux fonctions multiplicativement croissantes un résultat de Davenport et Erd\H os concernant les ensembles de multiples (1936, 1951, cf. \cite{zbMATH03023559}, \cite{zbMATH03065959}). Nos démonstrations suivent précisément l'approche élémentaire adoptée par ces auteurs dans \cite{zbMATH03065959} ; par souci d'autonomie, nous donnons néanmoins tous les détails.

\subsection{Décomposition en facteurs friable et criblé}

Soit $y>1$. Un nombre entier positif est dit $y$-friable\footnote{Dans certains textes en langue anglaise, la locution \emph{$y$-smooth} est utilisée pour désigner un tel nombre entier.} si tous ses diviseurs premiers sont $\ioe y$. L'ensemble de ces nombres sera noté $\Scal(y)$. Un nombre entier positif est dit $y$-criblé si tous ses diviseurs premiers sont $>y$. L'ensemble de ces nombres sera noté $\Ecal(y)$. Ces deux ensembles sont des facteurs directs complémentaires : tout nombre entier positif $n$ a une unique décomposition 
\begin{equation*}
n=ab, \text{ où } a \in \Scal(y) \text{ et }b \in \Ecal(y).
\end{equation*}
Nous dirons que $a$ est le facteur $y$-friable de $n$, et $b$ son facteur $y$-criblé.

Dans ce cas, le théorème d'Erd\H os, Saffari et Vaughan est immédiat ; on a les relations 
\begin{align}
\sum_{a \in \Scal(y)} \frac 1a&= \prod_{p \ioe y} (1-1/p)^{-1}\label{t19}\\
t^{-1} \sum_{\substack{b\ioe t\\b\in \Ecal(y}}1 &\vers \prod_{p \ioe y} (1-1/p) \quad (t \vers \infty).\label{t10}
\end{align}

Si $z\soe y>1$, le facteur $y$-friable de $n$ divise son facteur $z$-friable, le facteur $z$-criblé de $n$ divise son facteur $y$-criblé. De plus, si
\begin{equation}\label{t12}
n=a'b', \text{ où } a' \in \Scal(z) \text{ et }b' \in \Ecal(z),
\end{equation}
alors
\begin{equation*}
a'=ab'', \text{ où } b''=(a',b) \in \Ecal(y).
\end{equation*}
En particulier, $n$ et $a'$ ont même facteur $y$-friable.

\smallskip

Si $\Scal(y,z)$ désigne l'ensemble des nombres entiers dont tous les diviseurs premiers $p$ vérifient l'encadrement $y<p\ioe z$, tout élément de $\Scal(z)$ est, de façon unique, le produit d'un élément de~$\Scal(y)$ et d'un élément de $\Scal(y,z)$. L'identité suivante est une variante de \eqref{t19} :
\begin{equation}\label{t16}
\sum_{n\in \Scal(y,z)}\frac 1n=\prod_{y<p\ioe z} (1-1/p)^{-1}.
\end{equation}

\subsection{Réduction friable d'une fonction arithmétique}

Si $f$ est une fonction arithmétique, et $y>1$, la \emph{réduite $y$-friable} de $f$ est la fonction arithmétique définie par
\begin{equation}\label{t11}
f(n;y)= f(a) \quad ( n\in \Nat^* ; \, n=ab, \,a \in \Scal(y), \, b \in \Ecal(y)),
\end{equation}
autrement dit, $f(n;y)=f\big(n;\Scal(y)\big)$, avec la notation \eqref{t34}.

\smallskip

Si $z\soe y>1$ et si $g$ est la réduite $y$-friable de $f$, alors, la réduite $z$-friable de $g$ est égale à $g$. En effet, avec les notations de \eqref{t12}, on a
\begin{equation}\label{t14}
g(n;z)=g(a')=f(a';y)=f(n;y)=g(n),
\end{equation}
puisque $n$ et $a'$ ont le même facteur $y$-friable.

\smallskip

Explicitons le cas particulier de la proposition \ref{t31} obtenu en prenant $A=\Scal(y)$ et $B=\Scal(y)$.

\begin{prop}\label{t15}
Soit $f$ une fonction arithmétique à valeurs réelles, et soit~$y>1$. On suppose que
\begin{equation}\label{t13}
\sum_{a \in \Scal(y)} \frac{\min\big(f(a),0\big)}a > -\infty.
\end{equation}

Alors
\[
t^{-1}\sum_{n\ioe t} f(n;y) \vers \alpha(f;y) \quad ( t \vers \infty),
\]
où
\[
\alpha(f;y)=\prod_{p\ioe y} \Big(1-\frac 1p\Big)\sum_{a \in \Scal(y)}\frac{f(a)}a \quad (\in\, ]-\infty,\infty]).
\]

On a donc aussi
\begin{equation}\label{t17}
(\ln t)^{-1}\sum_{n\ioe t} \frac{f(n;y)}n \vers \alpha(f;y) \quad ( t \vers \infty).
\end{equation}
\end{prop}

\begin{prop}\label{t18}
Soit $f$ une fonction arithmétique à valeurs réelles, et soit~$z \soe y>1$. On suppose que \eqref{t13} est vérifiée et on considère la fonction arithmétique $g(n)=f(n;y)$ définie par \eqref{t11}. 

Alors, \eqref{t13} est vérifiée en remplaçant $f$ par $g$ et $y$ par $z$, et
\[
\alpha(g;z)=\alpha(f;y).
\]
\end{prop}
\dem

Compte tenu de \eqref{t16}, \eqref{t14} et \eqref{t13}, on a 
\begin{equation*}
\sum_{a' \in \Scal(z)} \frac{\min\big(g(a'),0\big)}{a' }=\sum_{\substack{a \in \Scal(y)\\b\in\Scal(y,z)}} \frac{\min\big(f(a),0\big)}{ab}=
\prod_{y<p \ioe z}(1-1/p)^{-1}\sum_{a \in \Scal(y)} \frac{\min\big(f(a),0\big)}a >-\infty.
\end{equation*}

Par conséquent,
\begin{equation*}
\alpha(g;z)=\lim_{t\vers\infty}t^{-1}\sum_{n\ioe t} g(n;z)=\lim_{t\vers\infty}t^{-1}\sum_{n\ioe t} f(n;y)=\alpha(f;y).\fine
\end{equation*}

\subsection{Existence de la valeur moyenne logarithmique des fonctions multiplicativement monotones}\label{t26}

Comme nous l'avons signalé ci-dessus, la proposition suivante généralise aux fonctions multiplicativement croissantes un théorème de Davenport et Erd\H os concernant les ensembles de multiples.

\begin{prop}\label{t8}
Soit $f$ une fonction arithmétique multiplicativement croissante. La quantité 
\[
\alpha(f;y)=\prod_{p\ioe y} \Big(1-\frac 1p\Big)\sum_{n \in \Scal(y)}\frac{f(n)}n 
\]
est une fonction croissante de $y$. Posons 
\[
\alpha =\lim_{y \vers \infty}\alpha(f;y) \quad (\in \, ]-\infty,+\infty]).
\]

Alors, on a
\[
\alpha=\lim_{x \vers \infty} \frac{1}{\ln x} \sum_{n\ioe x} \frac{f(n)}n=\liminf_{x \vers \infty} \frac 1x\sum_{n\ioe x}f(n).
\]
\end{prop}
\dem

Si $f$ est multiplicativement croissante, elle est bornée inférieurement, donc la condition \eqref{t13} est vérifiée et $\alpha(f;y)$ est bien défini pour tout $y>1$. C'est une fonction croissante de $y$. En effet, soit $z\soe y>1$. Comme $f$ est multiplicativement croissante, on a
\begin{equation*}
f(n;z) \soe f(n;y) \quad (n \in \Nat^*),
\end{equation*}
d'où $\alpha(f;z) \soe \alpha(f;y)$, d'après la proposition \ref{t15}.

\smallskip

Posons maintenant
\begin{align*}
\alpha &=\lim_{y\vers \infty}\alpha(f;y)=\sup_{y >1}\alpha(f;y)\\
\beta&=\liminf_{x\vers \infty}\frac 1x\sum_{n \ioe x} f(n)\\
\gamma &=\liminf_{x\vers \infty}\frac 1{\ln x}\sum_{n \ioe x} \frac{f(n)}n\\
\delta &=\limsup_{x\vers \infty}\frac 1{\ln x}\sum_{n \ioe x} \frac{f(n)}n\cdotp
\end{align*}

Les quatre limites $\alpha, \beta, \gamma, \delta$ appartiennent à $]-\infty,\infty]$. Nous allons voir qu'elles sont toutes égales en démontrant les inégalités
\[
\alpha \ioe \beta \ioe \gamma \ioe \delta \ioe \alpha.
\]
Parmi elles, l'inégalité $\gamma \ioe \delta$ est évidente, et l'inégalité $\beta \ioe \gamma$ est une relation, vraie pour toute fonction arithmétique $f$, qui se démontre par sommation partielle. Il reste à prouver les deux inégalités concernant $\alpha$. 

\smallskip

\underline{Preuve de l'inégalité $\alpha \ioe \beta$.}

Soit $y>1$. Pour tout $n \in \Nat^*$ on a $f(n) \soe f(n;y)$, donc 
\[
\beta=\liminf_{t\vers \infty}t^{-1}\sum_{n \ioe t} f(n) \soe \lim_{t\vers \infty}t^{-1}\sum_{n \ioe t} f(n;y)=\alpha(f;y),
\]
d'après la proposition \ref{t15}. Le résultat s'obtient en faisant tendre $y$ vers l'infini.

\smallskip

\underline{Preuve de l'inégalité $\delta \ioe \alpha$.}

Si $\alpha=\infty$, l'inégalité est évidente. Supposons donc $\alpha <\infty$. 

Soit $y>1$. Posons
\begin{align*}
g(n) &=f(n;y) \\
h(n)&=f(n)-g(n).
\end{align*}

Comme $f$ est multiplicativement monotone, on a $h \soe 0$. Par conséquent, pour $x>y$,
\begin{align}
\prod_{p\ioe x}(1-1/p)\sum_{n \ioe x} \frac{h(n)}n &\ioe\prod_{p\ioe x}(1-1/p) \sum_{n \in \Scal(x)} \frac{h(n)}n\notag\\
&=\alpha(h;x)=\alpha(f;x)-\alpha(g;x)\notag\\
&=\alpha(f;x)-\alpha(f;y)\expli{d'après la proposition \ref{t18}}\notag\\
&\ioe \alpha-\alpha(f;y).\label{t6}
\end{align}

Maintenant, quand $x$ tend vers l'infini,
\begin{align*}
\frac{1}{\ln x} \sum_{n\ioe x} \frac{f(n)}n &= \frac{1}{\ln x} \sum_{n\ioe x} \frac{g(n)}n + \frac{1}{\ln x} \sum_{n\ioe x} \frac{h(n)}n \\
&\ioe \alpha(f;y)+o(1)+\frac{1}{\ln x} \prod_{p\ioe x}(1-1/p)^{-1}\big(\alpha-\alpha(f;y)\big)\expli{d'après \eqref{t17} et \eqref{t6}}\\
&\ioe \alpha +o(1)+O\big(\alpha-\alpha(f;y)\big),
\end{align*}
où l'on a utilisé la majoration élémentaire $\prod_{p\ioe x}(1-1/p)^{-1} \ll \ln x$. On a donc, pour tout $y$, la majoration $\delta \ioe \alpha+O\big(\alpha-\alpha(f;y)\big)$, ce qui entraîne la conclusion en faisant tendre $y$ vers l'infini.~\fin

\smallskip

Si $f$ est multiplicativement monotone, nous poserons 
\[
\alpha(f)=\lim_{x \vers \infty}\frac{1}{\ln x} \sum_{n\ioe x} \frac{f(n)}n\cdotp
\]
C'est la \emph{valeur moyenne logarithmique} de $f$. Notons que $\alpha(f) \soe f(1)$, si $f$ est multiplicativement croissante, et $\alpha(f) \ioe f(1)$, si $f$ est multiplicativement décroissante.

\smallskip

L'égalité (8),~p.~22, de l'article \cite{zbMATH03065959} de  Davenport et Erd\H os, est généralisée aux fonctions multiplicativement croissantes dans la proposition suivante.
\begin{prop}
Soit $(f_k)_{k \soe 1}$ une suite croissante de fonctions arithmétiques multiplicativement croissantes. On suppose que 
\[
f(n)=\sup_{k \soe 1} f_k(n) < \infty \quad (n \in \Nat^*).
\]

Alors $\alpha(f_k)$ tend vers $\alpha(f)$ quand $k$ tend vers l'infini.
\end{prop}
\dem

On sait que $f$ est multiplicativement croissante. Pour tout $y>1$, $\alpha(f_k;y)$ tend vers $\alpha(f;y)$, quand $k$ tend vers l'infini, par convergence monotone. Par suite,
\begin{align*}
\lim_{k \vers \infty} \alpha(f_k) &=\sup_{k \soe 1} \alpha(f_k)=\sup_{k \soe 1}\,\sup_{y>1}\alpha(f_k;y)\\
&=\sup_{y>1}\,\sup_{k \soe 1}\alpha(f_k;y)=\sup_{y>1}\alpha(f;y)=\alpha(f).\fine
\end{align*}

\section{Remarques sur la moyenne de Cesàro des fonctions multiplicativement monotones}\label{t24}

En vertu de la proposition \ref{t8}, l'existence de la valeur moyenne d'une fonction arithmétique multiplicativement croissante, \cad de la limite, éventuellement infinie, de ses moyennes de Cesàro, équivaut à l'inégalité
\begin{equation}\label{t25}
\limsup_{x \vers \infty} \frac 1x\sum_{n\ioe x}f(n) \ioe \alpha(f).
\end{equation}

On sait depuis 1934 que cette inégalité n'est pas toujours vraie, puisque Besicovitch a montré dans \cite{zbMATH03014226} comment construire des ensembles de multiples sans densité asymptotique.

\medskip

Erd\H os a donné en 1948 (cf. \cite{zbMATH03049250}) un critère nécessaire et suffisant pour la validité de \eqref{t25} dans le cas d'un ensemble de multiples. Il serait intéressant de généraliser ce critère aux fonctions arithmétiques multiplicativement croissantes.

\medskip

Observons que l'existence de la valeur moyenne est évidente dans le cas où la dérivée au sens de Bougaïef est de signe constant :
\[
\frac 1x\sum_{n\ioe x}f(n)=\frac 1x\sum_{m\ioe x}Df(m)\lfloor x/m\rfloor \vers \sum_{m\soe 1}\frac{Df(m)}m \quad ( x \vers \infty).
\]

\section{Déterminants Toeplitz-multiplicatifs hermitiens positifs}\label{t3}

Toute matrice hermitienne définie positive infinie\footnote{J'entends par là que toutes les sections finies $(h_{ij})_{1\ioe i, j \ioe n}$ ($n \in \Nat^*$) sont définies positives.}
\[
H=(h_{ij})_{i, j \in \Nat^*}
\]
est la matrice de Gram 
\[
G=(\scal{e_i}{e_j})_{i,j \in \Nat^*}
\]
d'une suite $(e_n)_{n \in \Nat^*}$ de vecteurs d'un espace de Hilbert $\Hcal$, par exemple $\ell^2(\Nat^*)$. Cela se démontre en appliquant l'algorithme de Cholesky à $H$. 

Si $H$ est une matrice de Toeplitz (au sens usuel, additif), on a $h_{ij}=c_0(i-j)$, donc
\begin{equation}\label{t27}
\scal{e_i}{e_j} =\scal{e_{i+k}}{e_{j+k}}\quad (i,j \in \Nat^*, \, k \in \Nat).
\end{equation}

Si $H$ est une matrice Toeplitz-multiplicative, on a $h_{ij}=c(i/j)$, donc
\begin{equation}\label{t28}
\scal{e_i}{e_j} =\scal{e_{ki}}{e_{kj}}\quad (i,j,k \in \Nat^*).
\end{equation}

Notons
\[
D_{n}=\det(h_{ij})_{1 \ioe i, j \ioe n} \quad (n \in \Nat^*)
\]
avec $D_0=1$ par convention. Une formule classique de géométrie hilbertienne nous donne
\[
\frac{D_{n}}{D_{n-1}}=\frac{\Gram(e_1,\dots,e_n)}{\Gram(e_1,\dots,e_{n-1})}=\dist_{\Hcal}^2\big(e_n,V_{n-1}\big),
\]
où
\[
V_k=\Vect(e_1, \dots, e_k) \quad (k \in \Nat),
\]
(on a donc $V_0={0}$).

Si $H$ est une matrice de Toeplitz au sens usuel, on a, pour $n \soe 2$,
\begin{align*}
\frac{D_{n}}{D_{n-1}} &=\dist_{\Hcal}^2\big(e_n,V_{n-1}\big)\\
&\ioe \dist_{\Hcal}^2\big(e_n,\Vect(e_2, \dots, e_{n-1})\big)\\
&=\frac{\Gram(e_2,\dots,e_n)}{\Gram(e_2,\dots,e_{n-1})}\\
&=\frac{\Gram(e_1,\dots,e_{n-1})}{\Gram(e_1,\dots,e_{n-2})} \expli{d'après \eqref{t27}}\\
&=\frac{D_{n-1}}{D_{n-2}} \cdotp
\end{align*}

C'est le résultat de Fekete mentionné au \S\ref{t22}. La convergence de la suite $D_{n}/D_{n-1}$, qui en résulte, entraîne celle de la suite~$D_n^{1/n}$, vers la même limite. Si les $c_k$ sont les coefficients de Fourier d'une fonction $f$, de période $1$, continue et positive, Szeg\H o démontra dans le même article \cite{zbMATH02616579} une conjecture de P\'olya, affirmant que cette limite est la moyenne géométrique des valeurs de $f$ :
\[
D_n^{1/n} \vers \exp\Big(\int_0^1\ln f(t) \, dt\Big) \quad (n \vers \infty).
\]

Ce résultat fut le premier d'une nouvelle théorie, dont l'ouvrage récent \cite{nikolski2017} est un exposé actualisé.

\smallskip

Le raisonnement ci-dessus s'adapte de la façon suivante au cas où $H$ est une matrice Toeplitz-multiplicative (hermitienne, définie positive). Si $n=km$, on a
\begin{align*}
\frac{D_{n}}{D_{n-1}} &=\dist_{\Hcal}^2\big(e_{km},V_{km-1}\big)\\
&\ioe \dist_{\Hcal}^2\big(e_{km},\Vect(e_m,e_{2m}, \dots, e_{(k-1)m})\big)\\
&=\frac{\Gram(e_m,e_{2m},\dots,e_{(k-1)m},e_{km})}{\Gram(e_m,e_{2m},\dots,e_{(k-1)m})}\\
&=\frac{\Gram(e_1,\dots,e_{k})}{\Gram(e_1,\dots,e_{k-1})} \expli{d'après \eqref{t28}}\\
&=\frac{D_{k}}{D_{k-1}} \cdotp
\end{align*}

La suite des quotients $D_n/D_{n-1}$ est donc multiplicativement décroissante. La suite de leurs logarithmes l'étant également, on obtient la conclusion suivante par application de la proposition~\ref{t8} à la fonction
\[
f(n) =-\ln (D_n/D_{n-1}) \quad (n \in \Nat^*).
\]

\begin{prop}\label{t29}
Avec les hypothèses et les notations de la question du \S\ref{t22}, la moyenne logarithmique
\[
\frac{1}{\ln N}\sum_{k\ioe N} \frac 1k\ln(D_k/D_{k-1})
\]
converge vers un élément $\alpha$ de $[-\infty, \ln c(1)]$ quand $N$ tend vers l'infini. De plus,
\[
\limsup_{N \vers \infty} D_N^{1/N}=e^{\alpha}.
\]
\end{prop}

Si $\alpha = -\infty$, on a donc une réponse positive à la question du \S\ref{t22} : la limite cherchée est nulle. Si~$\alpha > -\infty$, on peut récrire l'assertion sur la moyenne logarithmique sous la forme
\[
\prod_{k=1}^N (D_k/D_{k-1})^{1/k}=N^{\alpha+o(1)} \quad ( N \vers \infty).
\]

\medskip

Cela étant, en toute généralité, je ne sais pas répondre à la question initiale. La construction de Besicovitch peut-elle être adaptée au cas des déterminants Toeplitz-multiplicatifs, et fournir ainsi un exemple où $D_n^{1/n}$ n'a pas de limite?

\medskip

Les propositions suivantes, \ref{t30} et \ref{t36}, indiquent des cas simples où cette limite existe.

\begin{prop}\label{t30}
Soit $A$ un facteur direct, stable par multiplication, et dont la somme des inverses converge. On pose
\begin{equation*}
R =\{a/a', \, a\in A, \, a' \in A\}
\end{equation*}

Soit $c : \, \Rat^{+*} \vers \Com$ une fonction telle que
\begin{align*}
c(1/r) &=\overline{c(r)} \quad (r \in R)\\
c(r) &=0 \quad (r\notin R),
\end{align*}
et telle que, pour tout $n \in \Nat^*$, le déterminant
\[
D_n =\det\big(c(i/j)\big)_{1\ioe i,j \ioe n}
\]
soit strictement positif. 

Alors la suite $D_n^{1/n}$ est convergente.
\end{prop}

\dem

Soit $B$ le facteur direct complémentaire de $A$. En interprétant $D_n$ comme déterminant de Gram, et en utilisant le fait que chaque nombre entier~$j\in~\Nat^*$ s'écrit de façon unique $j=ab$, avec~$a \in A$,~$b\in B$, on a
\[
D_n=\Gram(e_{ab} \, ; \,  a \in A, \, b\in B, \, ab\ioe n),
\]
où nous utilisons aussi le fait qu'un déterminant de Gram est indépendant de l'ordre des vecteurs qui le composent.

Si $b,b' \in B$, et $a,a',\alpha,\alpha' \in A$, on a
\[
\frac{ab}{a'b'}=\frac{\alpha}{\alpha'}\implique \alpha' a\cdot b= \alpha a' \cdot b' \implique b=b',
\]
puisque $A$ et $B$ sont des facteurs directs complémentaires et puisque $A$ est stable par multiplication. Comme $c$ est à support dans $R$, cela montre que
\[
\scal{e_i}{e_j}=0 \quad (i=ab\, ; \, j=a'b' \, ; \,  a,a' \in A \, ; \,  b,b' \in B \, ; \, b\neq b').
\]
Autrement dit, les vecteurs $e_i$ et $e_j$ sont orthogonaux si les décompositions de $i$ et $j$ en produits d'un élément de $A$ par un élément de $B$ ne font pas intervenir le même élément de $B$.  Cela permet de décomposer $D_n$ en produit :
\begin{align*}
D_n&=\prod_{ b \in B, b\ioe n} \Gram(e_{ab}\, ; \,  a \in A, \, a\ioe n/b)\\
&=\prod_{ b \in B, b\ioe n} \Gram(e_{a}\, ; \, a \in A, \, a\ioe n/b),
\end{align*}
d'après \eqref{t28}.

En écrivant maintenant $n=\alpha\beta$, avec $\alpha\in A$, $\beta \in B$, on a 
\begin{equation*}
\Gram(e_{a}\, ; \, a \in A, \, a\ioe n/b)=\Gram(e_{a}\, ; \, a \in A, \, a\ioe (n-1)/b) \quad ( b\neq \beta).
\end{equation*}
Par conséquent
\begin{align}
\frac{D_n}{D_{n-1}} &=\frac{\Gram(e_{a}\, ; \, a \in A, \, a\ioe n/\beta)}{\Gram(e_{a}\, ; \, a \in A, \, a\ioe (n-1)/\beta)}\notag\\
&=\frac{\Gram(e_{a}\, ; \, a \in A, \, a\ioe \alpha)}{\Gram(e_{a}\, ; \, a \in A, \, a\ioe \alpha -1)}\label{t37}\\
&=\frac{D_{\alpha}}{D_{\alpha-1}}\virg\notag
\end{align}
puisque \eqref{t37} ne dépend pas de $\beta$.

On peut donc appliquer la proposition \ref{t31} à la fonction arithmétique 
\[
f(n) =-\ln (D_n/D_{n-1}) \quad (n \in \Nat^*),
\]
puisqu'elle vérifie $f(n)=f(n;A)$ et $f(n) \soe -\ln c(1)$ pour tout $n$. Le résultat s'en déduit.\fin

\medskip

\begin{prop}\label{t36}
Soit $\sigma : \Nat^* \vers \Com$ une fonction multiplicative, et $c: \Rat^{+*} \vers \Com$ la fonction définie par
\[
c(i/j)=\sigma\big(i/(i,j)\big)\cdot\overline{\sigma\big(j/(i,j)\big)} \quad (i,j \soe 1),
\]

On suppose que, pour tout $n \in \Nat^*$, le déterminant
\[
D_n =\det\big(c(i/j)\big)_{1\ioe i,j \ioe n}
\]
est strictement positif. 

Alors la suite $D_n^{1/n}$ est convergente.
\end{prop}
\dem

Au point $(ii)$ du Corollary 3.4, p. 274 de l'article \cite{zbMATH05082236}, Hilberdink donne une expression de~$D_n$ :
\[
D_n=\prod_p\prod_{k\soe 1} (\rho_{k+1}(p)/\rho_k(p))^{\lfloor n/p^k\rfloor},
\]
où, pour tout nombre premier $p$, 
\[
\rho_{k}(p)=\frac{\Delta_k(p)}{\Delta_{k-1}(p)} \quad (k \soe 1),
\]
où $\big(\Delta_k(p)\big)_{k\soe 1}$ est une suite de déterminants de Toeplitz usuels (additifs), hermitiens et positifs. 

Le théorème de Fekete montre que les quotients $\rho_{k+1}(p)/\rho_k(p)$ sont tous $\ioe 1$. Par conséquent,
\[
\frac 1n \ln D_n=\sum_p\sum_{k\soe 1} \ln \big(\rho_{k+1}(p)/\rho_k(p)\big)\frac{\lfloor n/p^k\rfloor}n \vers \sum_p\sum_{k\soe 1} \frac{\ln \big(\rho_{k+1}(p)/\rho_k(p)\big)}{p^k}\quad (n \vers \infty)
\]

Dans ce cas, la dérivée de $\ln(D_n/D_{n-1})$ au sens de Bougaïef est $\ioe 0$. \fin

\smallskip

Terminons par deux remarques sur cette proposition \ref{t36}. D'abord, le Theorem 3.2, p. 274 de \cite{zbMATH05082236} montre que l'hypothèse de stricte positivité des déterminants $D_n$ équivaut à la stricte positivité des déterminants $\Delta_k(p)$ (notés $\det T_k(\widehat{f_p})$ par Hilberdink). Ensuite, la proposition \ref{t36} s'applique en particulier au cas, également considéré par Hilberdink, où la fonction $\sigma$ est complètement multiplicative. Dans ce cas, les matrices $M_n$ sont définies positives si, et seulement si $\lvert \sigma(p)\rvert <1$ pour tout nombre premier $p$. Le Corollary 5.1, p. 283 de \cite{zbMATH05082236}, montre alors que
\[
D_n^{1/n} \vers \prod_p\big(1-\lvert \sigma(p)\rvert ^2\big)^{1/p} \quad (n \vers \infty).
\]

\bigskip

\begin{center}
  {\sc Remerciements}
\end{center}
\begin{quote}
{\footnotesize 
La recherche ayant conduit à ces résultats a été subventionnée par le Conseil Européen de la Recherche (ERC), dans le Septième Programme-Cadre de l'Union Européenne (FP7/2007-2013), accord de subvention ERC  n$^{\text{o}}$ 670239.  

L'article a été, en partie, élaboré au printemps 2017, lors d'un séjour au département de mathématiques de l'University College de Londres. Je remercie Andrew Granville et cette institution pour les conditions de travail idéales dont j'ai alors bénéficié.}
\end{quote}

\providecommand{\bysame}{\leavevmode ---\ }
\providecommand{\og}{``}
\providecommand{\fg}{''}
\providecommand{\smfandname}{et}
\providecommand{\smfedsname}{\'eds.}
\providecommand{\smfedname}{\'ed.}
\providecommand{\smfmastersthesisname}{M\'emoire}
\providecommand{\smfphdthesisname}{Th\`ese}

%\bibliographystyle{smfplain}
%\bibliography{journal}{}

\medskip

\footnotesize

\noindent BALAZARD, Michel\\
Aix Marseille Université, CNRS, Centrale Marseille, I2M UMR 7373\\
13453 Marseille\\
FRANCE\\
Adresse \'electronique : \texttt{balazard@math.cnrs.fr}

\end{document}